\documentclass{amsart}

\usepackage{amsthm}
\usepackage{amsmath}
\usepackage{amssymb}
\usepackage{amsfonts}
\usepackage{amscd}
\usepackage{eufrak}
\usepackage{textcomp}
\usepackage{enumerate}

\begin{document}

    \title{About the logarithm function over the matrices}
    \author{GERALD BOURGEOIS}
    \address{Departement de Mathematiques, Faculte de Luminy, 163 avenue de
Luminy, case 901, 13288 Marseille Cedex 09, France.}
    \email{bourgeoi@lumimath.univ-mrs.fr}

\begin{abstract}
 We prove the following results: let $x,y$ be $(n,n)$ complex matrices
 such that $x,y,xy$ have no eigenvalue in $]-\infty,0]$ and $log(xy)=log(x)+log(y)$. If $n=2$, or if $n\geq3$ and $x,y$ are
  simultaneously triangularizable, then $x,y$ commute. In both cases we reduce the problem to a result in complex analysis.\\

\end{abstract}

\maketitle

    \section{Introduction}

$\mathbb{Z}^{*}$ refers to the non-zero integers.\\
 Let $u$ be a complex number. Then $Re(u),Im(u)$ refer to the real and imaginary parts of $u$; if $u\notin]-\infty,0]$
 then $arg(u)\in]-\pi,\pi[$ refers to  its principal argument.
\subsection{Basic facts about the logarithm.}
 Let $x$ be a complex $(n,n)$ matrix which
hasn't any eigenvalue in $]-\infty,0]$. Then $log(x)$, the $x$-principal logarithm, is the $(n,n)$ matrix $a$ such that:\\
 $e^a=x$ and the eigenvalues of $a$ lie in the strip $\{z\in\mathbb{C}: Im(z)\in]-\pi,\pi[\}$.\\
  $log(x)$ always exists and is unique; moreover $log(x)$ may be written as a polynomial in $x$.\\
\indent Now we consider two matrices $x,y$ which have no
eigenvalue in $]-\infty,0]$:\\
$\bullet$ If $x,y$ commute then $x,y$ are simultaneously
triangularizable and we may associate pairwise their eigenvalues
$(\lambda_j),(\mu_j)$; if moreover
$\forall{j},|arg(\lambda_j)+arg(\mu_j)|<\pi$, then $log(xy)=log(x)+log(y)$.\\
$\bullet$ Conversely if $xy$ has no eigenvalue in $]-\infty,0]$
and $log(xy)=log(x)+log(y)$ then do $x,y$ commute ? We will prove
that it's true for $n=2$ (theorem 1) or, for all $n$, if $x,y$ are
simultaneously triangularizable (theorem 2). But if $n>2$,
then we don't know the answer in the general case. \\
 \subsection{Lemma 1.} Let $x,y$ be two complex $(n,n)$ matrices
such that $x,y$ haven't any eigenvalue in $]-\infty,0]$ and $log(x)log(y)=log(y)log(x)$.\\
Then $x,y$ commute.\\

\textbf{Proof.} The principal logarithm over
$\mathbb{C}\;\backslash\;]-\infty,0]$ is
 one to one; thus, using Hermite's interpolation formula, $x$ or $y$ may be written as a polynomial in
 $log(x)$ or $log(y)$.$\;\square$\\

    \section{Dimension 2}

\subsection{Principle of the
proof.}
 \noindent The proof is based on the two next propositions. The first one is a corollary
of a Morinaga and
Nono's result (\cite [p. 356]{1}); the second is a technical result using complex analysis.\\

\textbf{Proposition 1.}
Let $\mathcal{U}=\{u\in\mathbb{C}^{*}:e^{u}=1+u\}$.\\
 Let $a,b$ be two $(2,2)$ complex matrices
such that $e^{a+b}=e^{a}e^{b}$ and $ab\not=ba$; let
$spectrum(a)=\{\lambda_1,\lambda_2\},spectrum(b)=\{\mu_1,\mu_2\}$.\\
\indent Then one of the three following $item$ is fulfilled:\\
(1)  $\lambda_1-\lambda_2\in{2i\pi\mathbb{Z}^{*}}$ and
$\mu_1-\mu_2\in{2i\pi\mathbb{Z}^{*}}$.\\
(2)  One of the following complex numbers $\pm(\lambda_1-\lambda_2)$, $\pm(\mu_1-\mu_2)$ is in $\mathcal{U}$.\\
(3)  $a$ and $b$ are simultaneously similar to
 $\begin{pmatrix}\lambda&0\\0&\lambda+u\end{pmatrix}$ and
 $\begin{pmatrix}\mu+v&1\\0&\mu\end{pmatrix}$ with
 $\lambda,\mu\in\mathbb{C}$, $u,v\in\mathbb{C}^{*},u\not=v$ and
 $\dfrac{e^{u}-1}{u}=\dfrac{e^{v}-1}{v}\not=0$.\\

 \textbf{Proposition 2.} Let $u,v$ be two distinct, non zero
complex numbers such that
$\dfrac{e^{u}-1}{u}=\dfrac{e^{v}-1}{v}\not=0$,
$|Im(u)|<2\pi,|Im(v)|<2\pi$.\\
 Then necessarily $|Im(u)-Im(v)|\geq{2\pi}$.\\

\textbf{Proof.} Assume that we can choose these $u,v$ such that
$|Im(u)-Im(v)|<{2\pi}$. Let $\lambda=\dfrac{e^{u}-1}{u}$ and let
$f$ be the holomorphic function: $f(z)=e^z-\lambda{z}-1$.\\
\indent Now we show that there exists $a\in]0,2\pi[$ such that
$Im(u),Im(v)$ are in $]-a,2\pi-a[$ and $f$ hasn't any zero with
imaginary part $-a$ or $2\pi-a$.\\
If it's false then $f$ admits an infinity of zeros in the strip
$\{z:Im(z)\in{]-2\pi,2\pi[}\}$:\\
Case 1: we can extract a sequence of zeros $z_k$ such that
$Re(z_k)\rightarrow-\infty$; then $f(z_k)\sim-\lambda{z_k}$,
a contradiction.\\
Case 2: we can extract a sequence of zeros $z_k$ such that
$Re(z_k)\rightarrow+\infty$; then $f(z_k)\sim{e^{z_k}}$, a
contradiction.\\
 \indent Let $r$ be a big positive real such that
$Re(u),Re(v)$ are in $]-r,r[$ and $\triangle$ be the rectangle
$\{z:-r\leq{Re(z)}\leq{r},-a\leq{Im(z)}\leq{2\pi-a}\}$. The
oriented edge $\partial\triangle$ consists of four parts:
$h_1=\{x+i(2\pi-a):x$ from $r$ to $-r\}$, $v_1=\{-r+iy:y$ from
$2\pi-a$ to $-a\}$, $h_2=\{x-ia:x$ from $-r$ to $r\}$,
$v_2=\{r+iy:y$ from $-a$ to $2\pi-a\}$. $f$ admits in $\triangle$
at least three zeros: $0,u,v$. Thus
$\dfrac{1}{2i\pi}\int_{\partial\triangle}\dfrac{f'(z)}{f(z)}dz=I(f(\partial\triangle),0)\geq{3}$
where $I$ refers to the index function.\\
\indent $f(z+2i\pi)-f(z)=-\lambda{2i\pi}$; then
$\{f(h_1),f(h_2)\}$ is in tubular form; moreover $f(h_1)$ and
$f(h_2)$ are isometric to a parametric curve in the form
$\{(e^t-\sigma{t},\tau{t}):t\in[-r,r]\}$ where $\sigma,\tau$ are
real; we can choose $a$ such that $\sigma,\tau\in\mathbb{R}^*$;
thus for $j\in\{1,2\}$
$|\dfrac{1}{2i\pi}\int_{h_j}\dfrac{f'(z)}{f(z)}dz|\leq{1}-\rho$ with $\rho=\dfrac{1}{2\pi}arctan(|\dfrac{\tau}{\sigma}|)$.\\
\indent We can choose $r,a$ such that $f'(z)\not=0$ on
$\partial\triangle$; thus $f(h_1),f(h_2)$ intersect
perpendicularly $f(v_1)$ and $f(v_2)$. If $z\in{v_1}$ and
$r\rightarrow{+\infty}$ then $f(z)=-\lambda{z}-1+O(e^{-r})$;
$f(v_1)$ is close to a segment of fixed direction and length. If
$z\in{v_2}$ and $r\rightarrow{+\infty}$ then $f(z)=e^{z}+O(r)$;
$f(v_2)$ is close
to an anticlockwise circle of radius $e^r$ containing $0$.\\
\indent Therefore
   $\dfrac{1}{2i\pi}(\int_{h_1}\dfrac{f'(z)}{f(z)}dz+\int_{h_2}\dfrac{f'(z)}{f(z)}dz)\leq{2-2\rho}$,
$\dfrac{1}{2i\pi}(\int_{v_1}\dfrac{f'(z)}{f(z)}dz)\approx{0}$,\\
$\dfrac{1}{2i\pi}(\int_{v_2}\dfrac{f'(z)}{f(z)}dz)\approx{1}$;
what is contradictory with
$I(f(\partial\triangle),0)\geq{3}$.$\;\;\square$

 \subsection{ Theorem 1.} Let $x,y$ be two $(2,2)$
complex matrices such that $x,y,xy$ haven't any eigenvalue in
$]-\infty,0]$ and $log(xy)=log(x)+log(y)$. Then $x,y$
commute.\\

\textbf{Proof.} We assume that $xy\not=yx$.
$e^{log(x)}e^{log(y)}=e^{log(x)+log(y)}$; using lemma 1,
$log(x)log(y)\not=log(y)log(x)$; thus we may use Proposition 1;
it's wellknown that $u\in\mathcal{U}$ implies that $|Im(u)|>2\pi$;
then, according to the logarithm definition, $a=log(x)$ and
$b=log(y)$ satisfy $item$ (3). Moreover the conditions
$|Im(u)|<2\pi,|Im(v)|<2\pi,|Im(u)-Im(v)|<2\pi$ are necessarily
fulfilled. Proposition 2 proves that these conditions can't be all
satisfied.\;\;\;$\square$\\


    \section{Dimension $n$}

$I$ refers to the identity matrix of dimension $n-1$. Let $\phi$
be the holomorphic function:
$\phi:z\rightarrow\dfrac{e^z-1}{z},\phi(0)=1$.\\
\textbf{Remark 1.} We have shown in part 2 that if $u,v$ are
complex numbers such that
$|Im(u)|<2\pi,|Im(v)|<2\pi,|Im(u-v)|<2\pi$ and
$\phi(u)=\phi(v)$, then $u=v$.\\
We'll use the following to prove our second main result.

\subsection{Proposition 3.} Let
$a=\begin{pmatrix}a_0&u\\0&\alpha\end{pmatrix},b=\begin{pmatrix}b_0&v\\0&\beta\end{pmatrix}$
be two complex $(n,n)$ matrices where $\alpha,\beta$ are complex
numbers and $a_0,b_0$ are $(n-1,n-1)$ complex matrices which
commute; let
$spectrum(a_0-\alpha{I})=(\alpha_i)_{i\leq{n-1}},spectrum(b_0-\beta{I})=(\beta_i)_{i\leq{n-1}}$.
     If $e^{a+b}=e^a{e^b}$ and $ab\not=ba$ then one of the following
$item$ must be satisfied:\\
(4) $\exists{i}:\beta_i\not=0$ and
$\phi(\alpha_i+\beta_i)=\phi(\alpha_i)$.\\
 (5) $\exists{i}:\alpha_i\not=0,\beta_i=0$ and
$\phi(-\alpha_i)=1$.\\

 \textbf{Proof.}  We
may assume that $a_0,b_0$ are upper triangular. Let
$a_1=a_0-\alpha{I},b_1=b_0-\beta{I},w=a_1v-b_1u=[w_1,\cdots,w_{n-1}]^T
$. Thus $ab\not=ba$ iff $w\not=0$.\\
 Then
 $e^a=e^\alpha\begin{pmatrix}e^{a_1}&\phi(a_1)u\\0&1\end{pmatrix},e^b=e^\beta\begin{pmatrix}e^{b_1}&\phi(b_1)v\\0&1\end{pmatrix}$,\\
$e^{a+b}=e^{\alpha+\beta}\begin{pmatrix}e^{a_1+b_1}&\phi(a_1+b_1)(u+v)\\0&1\end{pmatrix}$;
therefore $e^{a+b}=e^a{e^b}$ iff\\
(6)\;\;  $(\phi(a_1+b_1)-\phi(a_1))u=(e^{a_1}\phi(b_1)-\phi(a_1+b_1))v$.\\
$(e^{a_1}\phi(b_1)-\phi(a_1+b_1))b_1=(\phi(a_1+b_1)-\phi(a_1))a_1$
and (6) imply that
$(e^{a_1}\phi(b_1)-\phi(a_1+b_1))b_1v=(\phi(a_1+b_1)-\phi(a_1))a_1v=(\phi(a_1+b_1)-\phi(a_1))b_1u$;
thus\\
(7) $(\phi(a_1+b_1)-\phi(a_1))w=0$.\\
 We have also
$e^b=e^{-a}e^{a+b}$; then we can prove by the same method that\\
(8) $(\phi(b_1)-\phi(-a_1))w=0$.\\
There exists $k$ such that $w_k\not=0$ and if $j>k$ then $w_j=0$.
Therefore (7),(8) imply that
$\phi(\alpha_k+\beta_k)=\phi(\alpha_k)$ and
$\phi(\beta_k)=\phi(-\alpha_k)$; we are done except if
$\alpha_k=\beta_k=0$.\\
\indent Now we assume that $\alpha_k=\beta_k=0$.
$\phi(a_1+b_1)-\phi(a_1)=\dfrac{1}{2}b_1(I+P(a_1,b_1)),e^{a_1}\phi(b_1)-\phi(a_1+b_1)=\dfrac{1}{2}a_1(I+P(a_1,b_1))$
where $P$ is an analytic function, defined on $\mathbb{C}^2$, which satisfies $P(0,0)=0$. (6) can be
rewritten as $(I+P(a_1,b_1))w=0$. Therefore $(1+P(0,0))w_k=0$, a
contradiction.\;\;\;$\square$\\

    \subsection{ Theorem 2.} Let $x,y$ be $(n,n)$ complex matrices
    such that $x,y,xy$ haven't any eigenvalue in
$]-\infty,0]$ and $log(xy)=log(x)+log(y)$. If moreover $x,y$ are
simultaneously triangularizable then $xy=yx$.\\

\textbf{Proof.} We assume that $x,y$ are upper-triangular and
$xy\not=yx$; we prove inductively the result for $n\geq2$.
$x=\begin{pmatrix}x_0&?\\0&\lambda\end{pmatrix}$,
$y=\begin{pmatrix}y_0&?\\0&\mu\end{pmatrix}$ where $x_0,y_0$ are
$(n-1,n-1)$ upper triangular matrices which haven't any eigenvalue
in $]-\infty,0]$ and
$\lambda,\mu\in\mathbb{C}\setminus]-\infty,0]$.  The matrices
$a=log(x),b=log(y)$ are polynomials in $x$ or $y$, thus they are
upper-triangular in form
$a=\begin{pmatrix}log(x_0)&?\\0&log(\lambda)\end{pmatrix}$,
$b=\begin{pmatrix}log(y_0)&?\\0&log(\mu)\end{pmatrix}$. Thus
$log(x_0y_0)=log(x_0)+log(y_0)$; according to the recurrence
hypothesis $x_0y_0=y_0x_0$ and then
$log(x_0)log(y_0)=log(y_0)log(x_0)$. Moreover $e^{a+b}=e^a{e^b}$
and, from lemma 1, $ab\not=ba$.\\
\indent Now we use Proposition 3 with
$\alpha=log(\lambda),\beta=log(\mu),a_0=log(x_0),b_0=log(y_0)$.
Here $\alpha_i,\beta_i,\alpha_i+\beta_i$ have imaginary parts in
$]-2\pi,2\pi[$ and according to Remark 1,
$item$ (4),(5) can't be satisfied. $\;\;\;\square$\\

\indent We conclude with an easy result.
\subsection{Proposition 4.}Let $x,y$ be two positive definite hermitian $(n,n)$ matrices so
that $log(xy)=log(x)+log(y)$. Then $xy=yx$.\\

\textbf{Proof.} $log(xy)$ exists because
$spectrum(xy)\subset{]0,\infty[}$; $a=log(x),b=log(y)$ are
hermitian matrices such that $e^{a+b}=e^{a}e^b$. Moreover
$e^{a+b}=(e^{a+b})^*=e^b{e^a}$ and $e^{a}e^b=e^{b}e^a$ or $xy=yx$.
\;\;$\square$\\

\textbf{Remark.} It's wellknown that if $a,b$ are bounded self
adjoint operators on a complex Hilbert space, then
$e^{a+b}=e^{a}e^b$ implies that $ab=ba$. ( $cf.$ \cite [Corollary
1]{2}).

    \section{Conclusion} When $n=2$, we know how to characterize
    the complex $(n,n)$ matrices $a,b$ such that $ab\not={ba}$ and $e^{a+b}=e^a{e^b}$; it allowed us to bring back our problem to a result
    of complex analysis. Unfortunately, if $n\geq{3}$, the classification of such matrices is unknown. For this reason we
    can't prove, in this last case, the hoped result without supplementary assumption.\\

\textbf{Acknowledgement.}  The author would like to thank F.
Nazarov for his participation in the proof of proposition 2.


\begin{thebibliography}{99}

        \bibitem[1]{1}
K. Morinaga and T.Nono. On the non-commutative solutions of the
exponential equation $e^x e^y = e^{x+y}$. J. Sci. Hiroshima univ.
(A)17, (1954), 345-358.

\bibitem[2]{2}
 C. Schmoeger. Remarks on commuting exponentials in Banach
 algebras. Proceedings of the American Mathematical Society.
 Volume 127, Number 5, (1999), pages 1337,1338.





 \end{thebibliography}
\end{document}